\documentclass{amsart}

\usepackage[dvips]{graphics}

\newtheorem{theorem}{Theorem}[section]

\theoremstyle{definition}

\begin{document}

\parskip=0.5\baselineskip
\baselineskip=1.44\baselineskip

\title
{Triple homomorphisms of C*-algebras}
\date{August 30, 2004; to appear in Southeast Asian Bulletin of Mathematics}
\author{Ngai-Ching Wong}

\address{
Department of Applied Mathematics, National Sun Yat-sen
  University, and National Center for Theoretical Sciences, Kaohsiung, 80424, Taiwan, R.O.C.}

\email{wong@math.nsysu.edu.tw}

\dedicatory{In memory of our beloved friend, Kosita Beidar (1951.10.21--2004.3.9)}

\thanks{Partially supported  by
Taiwan National Science Council}



\begin{abstract}
 In this note, we will discuss what kind of operators between C*-algebras preserves Jordan triple products
 $\{a,b,c\}= (ab^*c + cb^*a)/2$.
 These include especially isometries and disjointness preserving operators.
 \end{abstract}

\maketitle

\section{Introduction}

Recall that a Banach algebra $A$ is an algebra with a norm $\|\cdot\|$ such that
$\|ab\|\leq \|a\|\|b\|$, and every Cauchy sequence converges.
A complex Banach algebra
$A$ is a C*-algebra if there is an involution ${}^*$ defined on $A$ such that $\|a^*a\|=
\|a\|^2$.
A special example is $B(H)$, the algebra of all bounded linear operators on a (complex)
Hilbert space $H$.
 By the Gelfand-Naimark-Segel
Theorem, C*-algebras are exactly those norm closed *-subalgebras of $B(H)$.
An abelian C*-algebra $A$ can also be represented as the algebra $C_0(X)$ of
continuous functions on
a locally compact Hausdorff space $X$ vanishing at infinity.  $X$ is compact if
and only if $A$ is unital.

It is well known that the algebraic structure determines the geometric (norm) structure
of a C*-algebra $A$.  Indeed, the norm of a self-adjoint element $a$ of $A$ coincides with
the spectral radius of $a$, and the latter is a pure algebraic object.  In
general, the norm of an
arbitrary element $a$ of $A$ is equal to $\|a^*a\|^{1/2}$, and $a^*a$ is self-adjoint.
For an abelian C*-algebra $A=C_0(X)$, we note that the underlying space $X$ can be considered as
the maximal ideal space of $A$ consisting of nonzero complex homomorphisms
(= linear and multiplicative functionals) of $A$.  The topology of $X$ is the hull-kernel
topology, and thus be solely determined by the algebraic structure of $A$.

In this note, we will discuss how much the algebraic structure can
be recovered if we know the norm, or other, structure of a
C*-algebra. In particular, isometries and disjointness preserving
operators of C*-algebras preserve triple products $\{a,b,c\}=
(ab^*c + cb^*a)/2$.

The author is very grateful to our late friend, Kosita Beidar,
from whom he learned how to look at a seemingly pure analytic
problem from the point of view of an algebraist.

\section{The geometric structure determines the algebraic structure}

Suppose $T: A\longrightarrow B$ is an isometric linear embedding between C*-algebras.
That is, $\|Tx\|=\|x\|$ for all  $x$ in $A$.  We are interested in knowing what kind of
algebraic structure $T$ inherits from $A$ to its range, which is in general just a Banach
subspace of $B$.
We begin with two famous results.

\begin{theorem} \emph{(Banach and Stone; see, e.g., \cite{JW})}
Let $X$ and $Y$ be locally compact Hausdorff spaces.
Let $T: C_0(X) \longrightarrow C_0(Y)$ be a surjective linear isometry.
Then $T$ is a weighted composition operator
$$
Tf=h\cdot f\circ\varphi, \quad \forall f\in C_0(X),
$$
where $h$ is a continuous scalar function on $Y$ with $|h(y)|\equiv 1$,
 and $\varphi$ is a homeomorphism from $Y$ onto $X$.
 Consequently, two abelian C*-algebras  are isomorphic as Banach spaces if, and only if,
they are isomorphic as ${}^*$-algebras.
\end{theorem}

Here is a sketch of the proof. Let $T^* : M(Y) \longrightarrow
M(X)$ be the dual map of $T$, which is again a surjective linear
isometry from the Banach space $M(Y)= C_0(Y)^*$ of all bounded
Radon measures on $Y$ onto that on $X$. Restricting $T^*$ to the
dual unit balls, which are weak* compact and convex, we get an
affine homeomorphism.  Since the nonzero extreme points of the dual unit
balls are exactly unimodular scalar multiples  of point masses,
$T^*$ sends a point mass $\delta_y$ to
$\lambda\delta_x$. Here $y\in Y$, $x\in X$ and $|\lambda|=1$.  We
write $x=\varphi(y)$ and $\lambda = h(y)$ to indicate that $x$ and
$\lambda$ depend on $y$.  It follows that
$$
Tf(y)=T^*(\delta_y)(f)=h(y)\delta_{\varphi(y)}(f)=h(y)f(\varphi(y)).
$$
It is then routine to see that $h$ is unimodular and continuous on $Y$, and that $\varphi$ is a homeomorphism
from $Y$ onto $X$.

\begin{theorem}[{Kadison \cite{K}}]
Let $A$ and $B$ be C*-algebras.  Let $T: A\longrightarrow B$ be a surjective linear isometry.
Then there is a unitary element $u$ in $\tilde B = B \oplus {\mathbb C}1$, the unitalization of $B$,
and a Jordan ${}^*$-isomorphism $J : A\longrightarrow B$ such that
$$
Ta = uJ(a), \quad\forall a\in A.
$$
Consequently, two C*-algebras are isomorphic as Banach spaces if, and only if, they are isomorphic
as Jordan ${}^*$-algebras.
\end{theorem}

Recall that a Jordan ${}^*$-homomorphism $J$ preserves linear sums, involutions and Jordan products:
$a\circ b = (ab + ba)/2$. It is easy to see that the abelian case can also be written in this form with
$u=h$ and $Jf=f\circ\varphi$.
In general, the product of a pair of elements in $A$ can be decomposed into two parts
$ab = a\circ b + [a,b]$, the sum of the Jordan product and the Lie product $[a,b]=(ab-ba)/{2}$.
It is plain that $a\circ b=b\circ a$ is commutative and $[a,b]=-[b,a]$ is anti-commutative.
However they are not associative.  The Kadison theorem states that
 the norm structure of a C*-algebra determines completely
its Jordan structure.

It is interesting to note that Jordan products are determined by squares:
$$
a\circ b = \frac{(a+b)^2 - a^2 - b^2}{2}, \quad \forall a,b\in A.
$$
A similar algebraic structure exists in C*-algebras, namely, the Jordan triple products:
$$
\{a,b,c\}= \frac{ab^*c+cb^*a}{2}.
$$
There is also a polar identity for triples:
$$
\{a,b,c\}= \frac{1}{8}\sum_{\alpha^2=1}\sum_{\beta^4=1}\, \alpha\beta\{a+\alpha b+\beta c\}^{(3)},
$$
Hence, a linear map $T$ between C*-algebras
 preserves triple products if and only if it preserves cubes
$a^{(3)}=\{a,a,a\}=aa^*a$.

Kaup \cite{KP} generalizes  Kadison theorem: a linear
surjection between C*-algebras  $T:A\longrightarrow B$ is an isometry
if and only if it preserves triple products.
A geometric proof of the Kadison Theorem is given by Dang, Friedman and Russo
\cite{DFR90}.
It goes first to note  that a norm exposed face of the dual unit ball
$U_{B^*}$ is of the form $F_u=\{\varphi\in B^*: \|\varphi\|=\varphi(u)\leq 1\}$ for
a unique partial isometry $u$ in $B^{**}$.  For two $\varphi, \psi$ in $B^*$, they
 are said to
be orthogonal to each other
if they have polar decompositions $\varphi=u|\varphi|, \psi = v|\psi|$ such
that $u\perp v$, i.e., $u^*v=uv^*=0$.  This amounts to say that $\|\varphi\pm\psi\|
=\|\varphi\|+\|\psi\|$.
Two faces $F_u, F_v$ are orthogonal if and only if $u\perp v$.
Then they verify that the dual map $T^*$ of  the surjective linear
isometry $T$  maps faces to faces
and preserves orthogonality. Consequently,
$T^{**} : A^{**}\longrightarrow B^{**}$ sends orthogonal partial isometries
to orthogonal partial isometries.
By the spectral theory, every element $a$ in $A\subset A^{**}$ can
be approximated in norm by a finite linear sum of orthogonal partial isometries
$\sum_j \lambda_j u_j$.  Then its cube $a^{(3)}$ can also be approximated by
$\sum_j \lambda_j^{(3)}u_j$.  It follows that $T(a^{(3)})$ and $(Ta)^{(3)}$ can
both be approximated by $\sum_j \lambda_j^{(3)}T^{**}u_j$.  Hence
$T(a^{(3)}) =(Ta)^{(3)}$, and thus $T$ preserves triple products by the polar identity.

We note that the above (geometric) proof of the Kadison theorem
depends very much on the fact the range of the
isometry is again a C*-algebra.
Extending the Holsztynski theorem
\cite{H, JW}, Chu and Wong \cite{CW} studied non-surjective linear isometries between
C*-algebras.

\begin{theorem}[{Chu and Wong \cite{CW}}]
Let $A$ and $B$ be C*-algebras and let $T$ be a linear isometry
from $A$ \emph{into} $B$.   There
is a largest closed projection $p$ in $B^{**}$ such that $T(\cdot)p : A
\longrightarrow B^{**}$ is a Jordan triple homomorphism and
$$
T(ab^*c + cb^*a)p= T(a)T(b)^*T(c)p + T(c)T(b)^*T(a)p, \quad \forall a,b,c\in A.
 $$
  When $A$ is abelian, we have
$\|T(a)p\|=\|a\|$ for all $a$ in $A$.
In particular,  $T$ reduces
{\it locally} to a Jordan triple isomorphism on the JB*-triple generated by any $a$ in $A$,
also an abelian C*-algebra,
by a closed  projection $p_a$.
\end{theorem}

Besides the triple technique,
the proof of the above theorem makes use of the concept of representing elements
of a C*-algebra
as special sections of a continuous field of Hilbert spaces developed in \cite{Wong94}.
It is still geometric.

\section{Disjointness preserving operators are triple
homomorphisms}

In this section, we do not assume the operator
$T$ is isometric. Although the following statement might have been
known to experts, we provide a new and short proof here as we do
not find any in the literature.

\begin{theorem}
Let $T: A\longrightarrow B$ be a bounded linear map between  C*-algebras.
Then $T$ is a triple homomorphism if and only if
$T^{**}$ sends partial isometries of $A^{**}$ to partial isometries of $B^{**}$.
\end{theorem}
\begin{proof}
One direction is trivial.
Suppose $T^{**}$ sends partial isometries of $A^{**}$ to partial isometries of $B^{**}$.
Let $u, v$ be two partial isometries in $A$.   Observe that they are orthogonal
to each other, namely, $u^*v=uv^*=0$, if and only if
they have orthogonal initial spaces and orthogonal range spaces.
This amounts to saying
that $u+\lambda v$ is a partial isometry for all scalar $\lambda$ with $|\lambda|=1$.
Consequently, $T$ sends orthogonal partial isometries to orthogonal partial isometries.
For every $a$ in $A\subset A^{**}$, approximate $a$ in norm by a finite linear sum $\sum_n \lambda_n u_n$ of
orthogonal partial isometries in $A^{**}$.
Then its cube $a^{(3)}=aa^*a$ can also be approximated in norm
by $\sum_n \lambda_n^{(3)} u_n$.  It follows that $Ta$ and $T(a^{(3)})$
can be approximated in norm by $\sum_n \lambda_n T^{**}u_n$
and $\sum_n {\lambda_n}^{(3)} T^{**}u_n$, respectively.
This gives $T(a^{(3)})=(Ta)^{(3)}$, $\forall a\in A$.  By the polar identity, we see that
$T$ is a triple homomorphism.
\end{proof}

We say that a linear map $T: A\longrightarrow B$ between C*-algebras is \emph{disjointness
preserving} if
$$
a^*b = ab^* = 0 \quad\text{implies}\quad (Ta)^*(Tb)=(Ta)(Tb)^* =0, \quad \forall a, b\in A.
$$
Clearly, $T$ is disjointness preserving if and only if it preserves disjointness of partial
isometries.
It is clear that every triple homomorphism preserves disjointness.
Looking at the well-known result of Jarosz  \cite{jarosz,JW}
in the abelian case,
we see that not every disjointness preserving map is a triple homomorphism.
Indeed, let $T: C_0(X)\longrightarrow C_0(Y)$ be a bounded disjointness preserving
linear map between abelian C*-algebras.  Then there is a closed subset $Y_0$ of $Y$ on which
every $Tf$ vanishes.  On $Y_1=Y\setminus Y_0$ there is a bounded continuous function $h$ and
a continuous map $\varphi$ from $Y_1$ into $X$ such that
$Tf_{\mid Y_1}=h\cdot f\circ\varphi$ for all $f$ in $C_0(X)$.
Hence, $T$ is a triple homomorphism if and only if $T^{**}1$
is a partial isometry in $C_0(Y)^{**}$.
We end this note with a proof of this fact for the non-abelian case.

\begin{theorem}
Let $T: A\longrightarrow B$ be a bounded linear map between
C*-algebras.  Then $T$ is a triple homomorphism if and only if
$T$ is  disjointness preserving and $T^{**}1$ is a partial isometry
in $B^{**}$.
\end{theorem}
\begin{proof}
We verify the sufficiency only.
For simplicity of notations, denote again by $T$
the bidual map of $T$ from $A^{**}$ into $B^{**}$.

Let $a$ be a self-adjoint element of  $A$.
Identify the commutative C*-subalgebra
of $A$ generated by $1$ and $a$ with $C(X)$, where $X \subseteq [0,\|a\|]$
is the spectrum of $a$.  Let $0=\alpha_0 < \alpha_1 < \cdots < \alpha_{n-1}
< \alpha_n=\|a\|+1$ such that $X=\cup_k X_k$ is a partition of $X$ with
$X_k=X\cap[\alpha_{k-1},\alpha_k)\neq\emptyset$, and
pick an arbitrary point $x_{k}$ from  $X_{k}$ for each $k=1,2,\ldots,n$.
In particular,
\begin{equation*}
1=\sum_{k} 1_{X_{k}},
\end{equation*}
where $1_{X_{k}}$ is the characteristic function of the
set ${X_{k}}$.
For $1<j< k$, we can
find two sequences $\{f_m\}_m$ and $\{g_m\}_m$ in $C(X)\cap A$ such that
$f_mg_{m+p} =0$ for $m,p=0,1,\ldots$, $f_m \to 1_{X_{j}}$
and $g_m \to 1_{X_{k}}$ pointwisely on $X$.  By the
weak$^*$ continuity of $T$, we see that
$$
T(1_{X_{k}})T(f_m)^* = \lim_{p\to\infty}
T(g_{m+p})T(f_m)^*=0 \quad\hbox{for all $m=1,2,\ldots$}.
$$
Thus
$$
T(1_{X_{k}})T(1_{X_{j}})^* = \lim_{m\to\infty}
T(1_{X_{k}})T(f_m)^*=0.
$$
Similarly, we have
$$
T(1_{X_{k}})^*T(1_{X_{j}}) = 0.
$$
When $j=1$, we note that $1_{X_1}$ is an open projection in $A^{**}$.
Hence there is an increasing net $\{a_\lambda\}_\lambda$ in $A$ converges
to $1_{X_1}$ weakly (see \cite[Proposition 3.11.9]{Ped}, and also \cite{Brown88}).
Using an argument similar to above, we still get
$$
T(1_{X_{k}})T(1_{X_{1}})^* = T(1_{X_{k}})^*T(1_{X_{1}}) = 0.
$$
Consequently, for each  $j=1,2,\ldots,n$,
we have
\begin{equation*}
T(1)T(1_{X_{j}})^*T(1) =
\sum_{m,n}T(1_{X_{n}})T(1_{X_{j}})^*T(1_{X_{m}})
=(T(1_{X_{j}}))^{(3)}.
\end{equation*}
This gives
$$
\sum_n T(1_{X_{n}})= T1 = (T1)^{(3)} = \sum_n (T(1_{X_{n}}))^{(3)}.
$$
Multiplying the above identity on the left by $T(1_{X_{n}})^*$ and
$((T(1_{X_{n}}))^{(3)})^*$ respectively,
we see that
$$
(T(1_{X_{n}}) - ((T(1_{X_{n}}))^{(3)})^*(T(1_{X_{n}}) - ((T(1_{X_{n}}))^{(3)})=0.
$$
Hence
$T(1_{X_{n}})$ is a partial isometry for each $n$ and orthogonal to the others.
It follows that
\begin{align*}
(T(f))^{(3)} &= \lim \left(\sum_{n}
f(x_{n})T(1_{X_{n}})\right)^{(3)}
=\lim \sum_{n} f(x_{n})^{(3)}(T(1_{X_{n}}))^{(3)}\\
&=\lim \sum_{n}
f(x_{n})^{(3)}T(1_{X_{n}}) =T(f^{(3)}),
\end{align*}
for all $f$ in $C(X)$.
By the polar identity, $T$ preserves triple products in $C(X)$.
Let $p=T1^*T1$ and $q=T1T1^*$ be the initial and range projections of the
partial isometry $T1$, respectively.  We have
$$
2Tf = 2T\{1,1,f\}= qTf + Tfp, \quad \forall f \in C(X).
$$
It then follows that
$$
Tf = qTf = Tfp,
$$
for all $f$ in $C(X)$.
Thus for each self-adjoint element $a$ of $A$,
\begin{equation}\label{qtfp}
T(a^{(3)})=(Ta)^{(3)} \quad\text{and}\quad Ta=qTa=Tap.
\end{equation}

At this point, we assume further that $T1=p=q$ is a projection.  It then
follows from \eqref{qtfp} that
for every self-adjoint element $a$ of $A$,
$$
Ta = T(a^*)= T\{1,a,1\}= \{T1,a,T1\} = T1(Ta)^*T1 = (Ta)^*.
$$
Therefore, $T$
also preserves self-adjointness.
For any self-adjoint elements $a, b$ of $A$, by observing the
identities $T(a\pm b)^{(3)} = (T(a\pm b))^{(3)}$, we also see
that $T(aba)=TaTbTa$.  Altogether, it follows that
$T(a+ib)^{(3)}=(T(a+ib))^{(3)}$.  In other words,
$T$ is a triple homomorphism from $A$ into $B$.

In general, we consider the map $J=T1^*T$
from $A$ into $T1^*TA\subseteq B^{**}$.  It follows from \eqref{qtfp} that
$J$ is also disjointness preserving.  Since $J1=T1^*T1=p$ is a projection,
$J$ is a triple homomorphism.  Hence $J(a^{(3)})=(Ja)^{(3)}$, and
thus
$$
T1^*T(a^{(3)})=T1^*Ta(Ta)^*T1T1^*Ta = T1^*Ta(Ta)^*qTa = T1^*(Ta)^{(3)},
$$
for all $a$ in $A$ by \eqref{qtfp}.  Therefore,
$$
qT(a^{(3)})=T1T1^*T(a^{(3)})=T1T1^*(Ta)^{(3)}=q(Ta)^{(3)},
$$
and then by \eqref{qtfp} again,
$$
T(a^{(3)})=(Ta)^{(3)}, \quad\forall a\in A.
$$
This says $T$ is a triple homomorphism from $A$ into $B$.
\end{proof}

\end{document}